\documentclass[graybox]{svmult}

\usepackage{mathptmx}       
\usepackage{helvet}         
\usepackage{courier}        
\usepackage{type1cm}        
\usepackage{makeidx}         
\usepackage{graphicx}        
\usepackage{multicol}        
\usepackage[bottom]{footmisc}

\usepackage{amsmath}
\usepackage{amssymb}
\usepackage{amscd}
\usepackage{indentfirst}

\def\N{{\Bbb N}}
\def\Z{{\Bbb Z}}

\def\1{{\bf 1}}

\makeindex             

\def\R{{\Bbb R}}
\def\vu{{\bf u}}
\def\vv{{\bf v}}
\def\vw{{\bf w}}
\def\vU{{\bf U}}

\def\vZ{{\bf Z}}
\def\vxi{{\bf \xi}}
\def\vN{{\bf \nabla}}

\begin{document}

\title*{Real  variable methods in harmonic analysis and Navier--Stokes equations}
\titlerunning{Harmonic Analysis and Navier--Stokes}
\author{Pierre Gilles LEMARI\'E-RIEUSSET}
\authorrunning{PG. Lemari\'e-Rieusset}
\institute{P.G. Lemari\'e-Rieusset \at LaMME, Univ Evry, CNRS, Universit\'e Paris-Saclay, 91025, Evry, France.\\ \email{pierregilles.lemarierieusset@univ-evry.fr} }

\maketitle

\abstract{Real variable methods in harmonic analysis were developed throughout the works of E.M. Stein. They turn out to be a powerful tool for the study of non-linear PDEs. We illustrate this point by discussing various points of the modern theory of Navier--Stokes equations.}


\section*{Introduction}

Among the seven Millenium problems proposed by the Clay Mathematical Institute, I shall consider the question of existence and smoothness of solutions to the Navier--Stokes equations.    Let us first recall the question raised by the Clay Mathematical Institute. as it has been presented by Ch. Fefferman in his 2000 talk at the Coll\`ege de France \cite{Fef06}~:
  \begin{quote}
 Let $\vu^0$ be any smooth, divergence-free vector field in the Schwartz class$\mathcal{S}(\R^3)$  Do there exist smooth functions  $p(t,x)$, $u_i(t,x)=(u_1(t,x),u_2(t,x),u_3(t,x))$ on $ \R^3\times [0,\infty)$ that satisfy 
$$ \partial_t\vu+(\vu\cdot\vN)\vu=\Delta\vu -\vN p,$$
$$  \vN\cdot \vu=0,$$
$$\vu(0,.)=\vu^0,$$
$$u, p\in\mathcal{C}^\infty(\R^3\times(0,+\infty))$$
and
$$ \sup_{t>0} \|\vu(t,.)\|_2\leq \|\vu^0\|_2 ?$$
\end{quote}

Commenting on the Clay Millenium Prize on Navier--Stokes equations, L. Tartar writes in 2006 
  \cite{Tar06} :
  \begin{quote}
  Reading the text of the conjectures
to be solved for winning that particular prize leaves the impression that the
subject was not chosen by people interested in continuum mechanics, as the
selected questions have almost no physical content.  /..../   The problems seem to have been chosen in the hope
that they will be solved by specialists of harmonic analysis,  and it has given the
occasion to some of these specialists to help others in showing the techniques
that they use, as in a recent book by Pierre Gilles LEMARI\'E-RIEUSSET\footnote{The book is the one I published in 2002 \cite{Lem02}.} /..../.
  \end{quote}
  The question I'd like to discuss here is to which extent harmonic analysis is used or should be used to study the Clay question on Navier--Stokes equations? Or, more generally, to discuss the Navier--Stokes equations on the whole space $\R^3$ in various functional settings while using tools from real-variable methods in harmonic analysis.

  The very first point is, of course, to define correctly what is called here harmonic analysis.  While   classical harmonic analysis has been devoted in the XIXth century  to the spectral analysis of the Laplace operator and of the heat equation, I shall focus on the theory that has been developed in the second half of the XXth century, mainly in the works of E. M. Stein \cite{Ste93}. (For a short account of this history, see the recent paper of  G. B. Folland \cite{Fol17}). As a matter of fact, the two Fields medalists who play influential roles on the Clay problem on Navier--Stokes equations, namely Ch. Fefferman and T. Tao, are both former students of E. M. Stein. A good account of this harmonic analysis theory is to be found in the books by Grafakos \cite{Gra08, Gra09}.

 The work of E. M. Stein is well illustrated by the titles of two of his books :  \textit{Singular Integrals and Differentiability Properties of Functions} \cite{Ste71}  and \textit{Harmonic Analysis~: Real-Variable Methods, Orthogonality, and Oscillatory Integrals} \cite{Ste93}. A major topic was the extension of Littlewood--Paley theory from the disc to $\R^n$. This is closely related to the study of Sobolev spaces and of Besov spaces, a class of spaces he studied thoroughly in his book on singular integrals \cite{Ste71}.
 
Littlewood--Paley--Stein decomposition of distributions and Besov spaces turned to be a fundamental tool for the modern approach of the Navier--Stokes equations and are the center of many books devoted to harmonic analysis and Navier--Stokes equations, such as  M. Cannone's \textit{Harmonic Analysis Tools for Solving the Incompressible Navier--Stokes equations} \cite{Can04}, H. Bahouri, J.Y. Chemin and R. Danchin's \textit{Fourier analysis and nonlinear partial differential equations}  \cite{BCD11} or P.G. Lemari\'e-Rieusset's \textit{Recent developments in the Navier--Stokes problem} \cite{Lem02}.

 But we shall try to show that the interaction of harmonic analysis with Navier--Stokes equations  is broader than the scope of Littlewood--Paley decomposition, and that many other ideas of E.M. Stein can be useful for future works. We shall ,pay a few words on the Clay problem, but as well on some points of the Navier--Stokes theory in more general settings such as Kato's mild solutions in $L^3$ (existence and uniqueness), or Serrin criteria for weak-strong uniqueness or regularity of Leray weak solutions.
 
 \section{Fourier transform.}

 \subsection*{Fourier--Navier--Stokes equations.}
 Naturally, Fourier transform plays an important role in the study of our problem, as it is a differential problem with constant coefficients and defined on the whole space.  If we note $\mathcal{F}_x$ the spatial Fourier transform
 $$ \mathcal{F}_x(f)(t,x)=\int_{\R^3} f(t,x) \, e^{-ix\cdot \xi}\, dx,$$
 the Navier--Stokes equations are turned into
 $$  \partial_t\mathcal{F}_x\vu(t,\xi)+\sum_{j=1}^3 i\xi_j  \mathcal{F}_x(u_j\vu)(t,\xi)=-\nu \vert\xi\vert^2 \mathcal{F}_x\vu(t,\xi) -i \mathcal{F}_xp(t,\xi) \, \vxi$$
 $$ \vxi\cdot\mathcal{F}_x\vu(t,\xi)=0$$
 $$ \mathcal{F}_x\vu(0,\xi)=\mathcal{F}_x\vu^0(\xi)= \vU^0(\xi).$$
 Moreover, we have
 $$ \mathcal{F}_x(u_j\vu)(t,\xi) =\frac 1{(2\pi)^3} \int_{\R^3} \mathcal{F}_xu_j(t,\eta) \mathcal{F}_x\vu(t,\xi-\eta)\, d\eta$$
 and 
 $$ \vert \xi\vert^2 \mathcal{F}_xp(t,\xi)=-\sum_{j=1}^3\sum_{k=1}^3 \xi_j\xi_k \mathcal{F}_x(u_ju_k)(t,\xi) .$$
 This gives the following  simple system on the vector $\mathcal{F}_x\vu=(\mathcal{F}_xu_1,\mathcal{F}_xu_2,\mathcal{F}_xu_3)$
 \begin{equation*}\begin{split} \partial_t\mathcal{F}_x u_l(t,\xi)=&-\vert\xi\vert^2 \mathcal{F}_x u_l(t,\xi) \\-& \sum_{j=1}^3 \sum_{k=1}^3 \int_{\R^3}  \frac{ i\xi_j\xi_k}{ (2\pi)^3\vert \xi\vert^2}  (\xi_l \mathcal{F}_xu_k(t,\eta)-\xi_k \mathcal{F}_xu_l(t,\eta)) \mathcal{F}_xu_j(t,\xi- \eta) \, d\eta.
 \end{split}\end{equation*}
 This is turned into an integral equation on $\vU=\mathcal{F}_x \vu$ :
 $$ \vU(t,\xi)=e^{-t\vert\xi\vert^2}\vU^0(\xi)- B(\vU,\vU)(t,\xi)$$
 with
 $$ B(U,V)_l(t,\xi)=  \int_0^t e^{-(t-s)\vert\xi\vert^2}   \sum_{j=1}^3 \sum_{k=1}^3 \int_{\R^3}  \frac{ i\xi_j\xi_k}{ (2\pi)^3\vert \xi\vert^2}  (\xi_l  U_k(s,\eta)-\xi_k U_l(s,\eta)) V_j(s,\xi- \eta) \, d\eta\, ds.$$

 This allows very simple computations for the search of solutions. Indeed, let us assume that $\vU^0$ is controlled by a function $W^0$ :
 $$  \vert \vU^0(\xi)\vert\leq W^0(\xi)$$
 and that $W(t,\xi)$ is measurable, almost everywhere finite and is a non-negative  solution of the integral inequation  for every $t\in [0,T]$ and every $\xi\in \R^3$
 $$  e^{-t\vert\xi\vert^2} W^0(\xi) + B_0(W,W)(t,\xi) \leq W(t,\xi)$$
 with
 $$ B_0(W,V)(t,\xi)=\frac {18}{(2\pi)^3}\int_0^t e^{-(t-s)\vert\xi\vert^2}  \vert \xi\vert\int_{\R^3} W(s,\eta) V(s,\xi-\eta) \, d\eta\, ds.$$
 Define $W^{[0]}(t,\xi):= e^{-t\vert\xi\vert^2} W^0(\xi)$, $W^{[n+1]}(t,\xi):= W^{[0]}(t,\xi)+B_0(W^{[n]},W^{[n]})(t,\xi)$ and similarly
 $\vU^{[0]}:=e^{-t\vert\xi\vert^2} \vU^0(\xi)$ and  $\vU^{[n+1]}(t,\xi):= \vU^{[0]}(t,\xi)-B(\vU^{[n]},\vU^{[n]})(t,\xi)$. By induction on $n$, we find that we have the pointwise inequalities
 \begin{itemize}
 \item[$\bullet$]  $0\leq W^{[n]}(t,\xi)\leq W^{[n+1]}(t,\xi)\leq W(t,\xi)$
  \item[$\bullet$]  $\vert U^{[n]}(t,\xi)\vert\leq W^{[n]}(t,\xi)$
  \item[$\bullet$]  $\vert U^{[n+1]}(t,\xi)-U^{[n]}(t,\xi)\vert\leq W^{[n+1]}(t,\xi)-W^{[n]}(t,\xi)$.
\end{itemize}
 We find that $W^{[n]}$ is pointwise convergent to  a function $W^{[\infty]}\leq W$. By monotonous convergence, we have
 $$ W^{[\infty]}=W^{[0]} +B_0(W^{[\infty]},W^{[\infty]}).$$
 Then, by dominated convergence, we find that $\vU^{[n]}$ converges to a limit $\vU^{[\infty]}$ such that
 $$ \vU^{[\infty]}=\vU^{[0]} -B(\vU^{[\infty]},\vU^{[\infty]}).$$ $\vU^{[\infty]}$ is then the Fourier transform of a solution to the Navier--Stokes problem with initial value $\vu_0$.

 \subsection*{Gevrey analyticity.}
 This formalism allows one to get Gevrey-type analyticity estimates. If we assume more precisely that
 $$ \vert\vU^0(\xi)\vert \leq \frac 1 {2e}   W^0(\xi)$$  then we find that, for $0\leq t\leq T$,
\begin{equation}\label{Gevrey}\vert\vU^{[\infty]}(t,\xi) \vert\leq \frac 1 {2\sqrt e}  e^{-\sqrt t \vert\xi\vert} W^{[\infty]}(\frac 1 2  t,\xi).\end{equation}
 Indeed, we write 
 $$ \sup_{z\geq 0} e^{z-\frac 1 2 z^2}= \sqrt e$$
 and, for $0\leq s\leq t$
 $$ e^{\sqrt t \vert \xi\vert} e^{-\sqrt s\vert\xi-\eta\vert}e^{-\sqrt s \vert\eta\vert}\leq e^{(\sqrt t - \sqrt s) \vert\xi\vert}\leq e^{\sqrt{t-s} \vert\xi\vert}.$$
 We define
 $$ \vZ^{[n]} (t,\xi)=e^{\sqrt t \vert\xi\vert} \vU^{[n]}(t,\xi).$$
 We have
 \begin{equation*}\begin{split}
 Z^{[n+1]}  =Z^{[0]} -B_*(\vZ^{[n]},\vZ^{[n]})
 \end{split}
 \end{equation*}
 with, for $\vZ=e^{\sqrt t\vert\xi\vert}\vU$, 
  \begin{equation*}\begin{split}
 B_*(\vZ,\vZ)_l(t,\xi) &=e^{\sqrt t \xi\vert} \int_0^t e^{-(t-s)\vert\xi\vert^2}   \sum_{j=1}^3 \sum_{k=1}^3 \int_{\R^3}  \frac{ i\xi_j\xi_k}{ (2\pi)^3\vert \xi\vert^2}  (\xi_l  U_k(s,\eta)-\xi_k U_l(s,\eta)) U_j(s,\xi- \eta) \, d\eta\, ds\\ = \int_0^t e^{-(t-s)\vert\xi\vert^2}  & \sum_{j=1}^3 \sum_{k=1}^3 \int_{\R^3} 
  e^{\sqrt t\vert \xi\vert -\sqrt s \vert\xi-\eta\vert -\sqrt s \vert\eta\vert} 
  \frac{ i\xi_j\xi_k}{ (2\pi)^3\vert \xi\vert^2}  (\xi_l  Z_k(s,\eta)-\xi_k Z_l(s,\eta)) Z_j(s,\xi- \eta) \, d\eta\, ds \end{split}
 \end{equation*}
 If $\vert \vZ(t,\xi)\vert\leq A(t/2,\xi)$, we find
  \begin{equation*}\begin{split}
 \vert B_*(\vZ,\vZ)(t,\xi) \vert&\leq \sqrt e \frac {18}{(2\pi)^3} \int_0^t e^{-\frac 1 2 (t-s)\vert\xi\vert^2}  \vert \xi\vert\int_{\R^3} A(s/2,\eta) A(s/2,\xi-\eta) \, d\eta\, ds
 \\ &\leq 2 \sqrt e \frac {18}{(2\pi)^3} \int_0^{t/2} e^{-   (\frac t 2-\sigma)\vert\xi\vert^2}  \vert \xi\vert\int_{\R^3} A(\sigma,\eta) A(\sigma,\xi-\eta) \, d\eta\, d\sigma
 \\ & = 2 \sqrt e B_0(A,A)(\frac t 2,\xi).
 \end{split}
 \end{equation*}
 By induction on $n$, we then find that
 $$ \vert \vZ^{[n]}(t,\xi)\vert\leq \frac 1 {2\sqrt e}  W^{[n] }(\frac t 2,\xi).$$
 Thus, we have proved the Gevrey estimate (\ref{Gevrey}).
 
 \subsection*{Cheap Navier--Stokes equation.}
 
 Thus far, we have introduced, as a tool for studying the Navier--Stokes problem 
\begin{equation} \label{NSE} \left\{  \begin{split}\partial_t\vu+(\vu\cdot\vN)\vu=\Delta\vu -\vN p\\
  \vN\cdot \vu=0\\
 \vu(0,.)=\vu^0, \end{split}\right.\end{equation}
the study of the equation
\begin{equation} \label{Fcheap}  W=W^{[0]}+B_0(W,W) \end{equation}
 or, taking the inverse Fourier transform $w=\mathcal{F}_x^{-1}W$ of $W$,  the equation
 \begin{equation} \label{cheap}\left\{  \begin{split}\partial_t w =\Delta w+ 18 \sqrt{-\Delta}(w^2)\\
 w(0,.)=w^0. \end{split}\right.\end{equation}
 Equation (\ref{cheap}) is known as the \textit{cheap Navier--Stokes equation}. It has been introduced in 
 2001 by S. Montgomery-Smith \cite{MoS01} as a toy model for the Navier--Stokes equations.  He gave an example of an initial value $w^0$ in the Schwartz class ($w^0\in\mathcal{S}(\R^3)$) with a non-negative Fourier transform $W^0$ such that the solution $w$ blows up in finite time.
 
 The study of   equation (\ref{Fcheap}) has provided simple classes of solutions  to the Navier--Stokes equations.
  For instance, if $Z(\xi) $ is a non-negative measurable function that   satisfies the following inequation
 $$  W^0(\xi)+\frac {18}{(2\pi)^3 \vert\xi\vert} \int_{\R^3} Z(\xi-\eta)Z(\eta)\, d\eta\leq Z(\xi), $$ we get, by induction on $n$, that $W^{[n]}(t,\xi)\leq Z(\xi)$. This means that, if $W^0$ belongs to a  lattice Banach space of functions $E$ such that the operator $(Z,V)\mapsto \frac 1{\vert \xi\vert} (Z*V)$ is bounded on $E$ and if $\|W^0\|_E$ is small enough, then the Navier--Stokes equations (\ref{NSE}) with initial value $\vu^0$ with $\vert \mathcal{F}_x \vu^0\vert\leq W^0$ has a global solution $\vu$ with $\sup_{0<t<+\infty} \vert \mathcal{F}_x \vu\vert\in E$. Two simple instances can be found in the litterature :
 \begin{itemize}
 \item[$\bullet$] the case where $E=L^2( \vert\xi\vert \, d\xi)$ : if $Z\in E$, it means that $Z=\frac 1 {\vert\xi\vert^{1/2}} Z_0$ with $Z_0\in L^2$; thus $Z$ belongs to the Lorentz space $L^{3/2, 2}$ (as a product of a function in $L^{6,\infty}$ by a function in $L^2$), thus $Z*Z$ belongs to $L^{3,1}\subset L^{3,2}$ and $\frac 1{\vert \xi\vert^{1/2}} Z*Z\in L^{2,2}=L^2$. Thus, we find that if the initial value $\vu_0$ has a small norm in the homogeneous Sobolev space $\dot H^{1/2}=\mathcal{F}_x^{-1}(L^2(\vert\xi\vert\, d\xi))$, then the Navier--Stokes problem with initial value $\vu_0$ has a global solution. This is the result of Fujita and Kato \cite{FuK62}.
 \item[$\bullet$] the equality
 $$ \int \frac 1{\vert \xi-\eta\vert^2} \frac 1 {\vert \eta\vert^2}\, d\eta = C_0 \frac 1{\vert \xi\vert}$$ allowed Le Jan and Sznitman to consider the space $E$ defined by $$Z\in E\Leftrightarrow  Z\in L^1_{\rm loc}\text{ and }  \vert\xi\vert^2 Z\in L^\infty.$$ Again, they found that,   if the initial value $\vu_0$ has a small norm in the homogeneous Besov space $\dot B^{-2}_{PM,\infty}=\mathcal{F}_x^{-1}( \frac 1 {\vert\xi\vert^2} L^\infty (d\xi))$, then the Navier--Stokes problem with initial value $\vu_0$ has a global solution \cite{LJS97}.
   \end{itemize}
   
   If we look for local-in-time solutions, we must include the time variable in our estimations. For instance, since $$ e^{-(t-s)\vert\xi\vert^2} \leq e^{\frac 3 4} \left(\frac 3 4\right)^{3/2} \frac 1 {(t-s)^{3/4}} \frac 1 {\vert \xi\vert^{3/2}}, $$ then, 
  if $Z(\xi) $ and $\alpha(t)$  are   non-negative measurable functions that   satisfy  on $(0,T)$ the following inequation
 $$  W^{[0]}(\xi)+ e^{\frac 3 4} \left(\frac 3 4\right)^{3/2} \frac {18}{(2\pi)^3}\int_0^t \frac {\alpha(s)^2} {(t-s)^{3/4}} \, ds \, \frac 1{  \vert\xi\vert^{1/2}} \int_{\R^3} Z(\xi-\eta)Z(\eta)\, d\eta \leq \alpha(t) \,Z(\xi), $$ we get, by induction on $n$, that $W^{[n]}(t,\xi)\leq  \alpha(t)\, Z(\xi)$.
  Thus, if  $F$ is  a  lattice Banach space of functions  such that the operator $(Z,V)\mapsto \frac 1{\vert \xi\vert^{1/2}} (Z*V)$ is bounded on $F$ and if   $\sup_{0<t<T} t^{1/4} e^{-t\vert\xi\vert^2} W^0(\xi) \in F$   and  $\|\sup_{0<t<T} t^{1/4} e^{-t\vert\xi\vert^2} W^0(\xi)\|_F$ is small enough,
   then the Navier--Stokes equations (\ref{NSE}) with initial value $\vu^0$ with $\vert \mathcal{F}_x \vu^0\vert\leq W^0$ has a global solution $\vu$ with $\sup_{0<t<T}  t^{1/4}\vert \mathcal{F}_x \vu\vert\in F$.  Let us look at our two simple instances :
 \begin{itemize}
 \item[$\bullet$] the case where $E=L^2( \vert\xi\vert \, d\xi)$ and $F=L^2(\vert\xi\vert^2\, d\xi)$ : if $Z\in F$, it means that $Z=\frac 1 {\vert\xi\vert} Z_0$ with $Z_0\in L^2$; thus $Z$ belongs to the Lorentz space $L^{6/5, 2}$ (as a product of a function in $L^{3,\infty}$ by a function in $L^2$) and $V=\vert \xi\vert^{1/2}\in L^{3/2,2}$, thus, writing
 $$ \frac 1 {\vert\xi\vert^{1/2}} \vert Z*Z\vert \leq \frac  2  {\vert \xi\vert} (\vert Z\vert*\vert V\vert ),$$
 we get that 
  $ \vert Z\vert*\vert V\vert$ belongs to $L^{6/5,2}*L^{3/2,2}\subset L^{2,1}\subset L^{2,2}=L^2$, so that  we have $\frac 1{\vert\xi\vert^{1 /2}} (Z*Z)\in F$.  Moreover, if $A>0$ and $W_0\in E$, we find that, for $t>0$,
  $$\vert  t^{1/4} e^{-t\vert\xi\vert^2}  W_0(\xi) \vert \leq  1_{\vert \xi\vert\leq A} t^{\frac 1 4} W_0(\xi)+ 1_{\vert\xi\vert>A}  \frac 1 {\vert\xi\vert^{1/2} } W_0(\xi)$$ so that
   $$\|\sup_{0<t<T} t^{1/4} e^{-t\vert\xi\vert^2} W^0(\xi)\|_F\leq T^{1/4} A^{1/2} \|W^0\|_E+ \| 1_{\vert\xi\vert>A}  W^0\|_E$$
   and
   $$\lim_{T\rightarrow 0^+} \|\sup_{0<t<T} t^{1/4} e^{-t\vert\xi\vert^2} W^0(\xi)\|_F=0.$$
  
 Thus, we find that if the initial value $\vu^0$ belongs to the homogeneous Sobolev space $\dot H^{1/2}=\mathcal{F}_x^{-1}(L^2(\vert\xi\vert\, d\xi))$, then the Navier--Stokes problem with initial value $\vu^0$ has a local in time  solution. This is the result of Fujita and Kato \cite{FuK62}. Gevrey regularity estimates for a data in the Sobolev spave were first given by Foias and Temam \cite{FoT89}.
 \item[$\bullet$] the case where $E=\frac 1  {\vert\xi\vert^2} L^\infty(d\xi)$ and $F=\frac 1 {\vert\xi\vert^{5/2}} L^\infty(d\xi)$ :  the equality
 $$ \int \frac 1{\vert \xi-\eta\vert^{5/2}} \frac 1 {\vert \eta\vert^{5/2}}\, d\eta = C_0 \frac 1{\vert \xi\vert^2}$$ 
 shows that $ \frac 1{\vert \xi\vert^{1/2}}(F*F)\subset F$.  Moreover, if $A>0$ and $W_0\in E$, we write again that, for $t>0$,
  $$\vert  t^{1/4} e^{-t\vert\xi\vert^2}  W_0(\xi) \vert \leq  1_{\vert \xi\vert\leq A} t^{\frac 1 4} W_0(\xi)+ 1_{\vert\xi\vert>A}  \frac 1 {\vert\xi\vert^{1/2} } W_0(\xi)$$ so that
   $$\|\sup_{0<t<T} t^{1/4} e^{-t\vert\xi\vert^2} W^0(\xi)\|_F\leq T^{1/4} A^{1/2} \|W^0\|_E+ \| 1_{\vert\xi\vert>A}  W^0\|_E$$
   and
   $$\limsup_{T\rightarrow 0^+} \|\sup_{0<t<T} t^{1/4} e^{-t\vert\xi\vert^2} W^0(\xi)\|_F\leq \limsup_{A\rightarrow +\infty} \sup_{\vert\xi\vert>A} \vert \xi\vert^2 \vert W_0(\xi)\vert.$$ Again, we find that,   if the initial value $\vu_0$ belongs to the homogeneous Besov space $\dot B^{-2}_{PM,\infty} $ and if  $\limsup_{A\rightarrow +\infty} \sup_{\vert\xi\vert>A} \vert \xi\vert^2 \vert \mathcal{F}_x \vu^0(\xi) \vert.$ is small enough, then the Navier--Stokes problem with initial value $\vu_0$ has a local in time solution.
   \end{itemize}

$\ $\\

We have considered the basic examples of $E=L^2(\vert\xi\vert\, d\xi)$ or $E=\frac 1{\vert\xi\vert^2} L^\infty(d\xi)$. But many other examples are known. In particular, the theory has been developed for $W^0$ in certain Herz spaces. Recall that the Herz space $\mathcal{B}^s_{p,q}$  \cite{Her68} is defined by
$$ W\in \mathcal{B}^s_{p,q}\Leftrightarrow (2^{js} \|1_{2^j\leq \vert \xi\vert<2^{j+1}} W\|_p)_{j\in\Z} \in l^q.
$$
For instance, we have $L^2(\vert\xi\vert\, d\xi)=\mathcal{B}^{1/2}_{2,2}$ and $\frac 1 {\vert \xi\vert^2} L^\infty(d\xi)=\mathcal{B}^2_{\infty,\infty}$.   In 2012, Cannone and Wu  \cite{CaW12} have studied the Navier--Stokes problem with an initial value $\vu_0$ such that $\mathcal{F}_x\vu_0\in \mathcal{B}^{-1}_{1,q}$ with $1\leq q\leq 2$. The case $q=1$ corresponds to the case $\mathcal{F}_x \vu_0\in L^1(\frac{d\xi}{\vert\xi\vert})$, a case studied by Lei and Lin in 2011 \cite{LeL11}.

\section{Singular integrals.}
\subsection*{Helmholtz decomposition.}
Modern history of  harmonic analysis begins with the study of singular integrals, from the work of M. Riesz  on the Hilbert transform in 1924 \cite{Rie24}  to the fundamental paper of Calder\'on and Zygmund on singular integrals in 1952 \cite{CaZ52} (and its extension to vector-valued integrals by Benedek, Calder\'on and Panzone in 1962 \cite{BCP62}). Basic accounts of the theory are to be found in the first chapters of the books of Stein \cite{Ste71, Ste93} or Grafakos \cite{Gra08}.

  The paradigm of Calder\'on--Zygmund convolution operators on $\R^d$ is given by   Marcinkiewicz multipliers: if $K$ is the inverse Fourier transform of a function $m(\xi)$ such that, for every $\alpha\in\N_0^d$ with $\vert \alpha\vert\leq d+2$, $$ \sup_{\xi\neq 0} \left\vert  \vert\xi\vert^{\vert \alpha \vert}\frac{\partial^\alpha m}{\partial\xi^\alpha}(\xi) \right\vert<\infty,$$ 
  then convolution with $K$ is a Calder\'on--Zygmund  operator. The most classical example is given by
   the Riesz transforms $R_j$, $j=1,\dots,d$ :
  $$ R_j f= \frac{\partial_j}{\sqrt{-\Delta}} f=\mathcal{F}_x^{-1} \left( \frac{i\xi_j}{\vert\xi\vert} \mathcal{F}_xf\right).$$
  Riesz transforms are naturally encountered when studying the Helmholtz decomposition of a vector field defined on the whole space $\R^3$. One considers a vector field $\vu$ and we want to decompose it as a sum of a divergence-free vector field $\vv$   and an irrotational vector field $\vw$   :
  $$ \vu=\vv+\vw \text{ with } \vN\cdot\vv \text{ and } \vN\wedge\vw=0.$$
 Basic formulas of vector analysis link the divergence and the curl of a vector $\vu$ to its Laplacian by
 $$ \vN\wedge (\vN\wedge \vu)= -\Delta \vu +\vN(\vN\cdot\vu). $$   
 In particular, we have
 $$ -\Delta\vv =\vN\wedge(\vN\wedge\vv)= \vN\wedge(\vN\wedge\vu)$$
 and
 $$ -\Delta w= -\vN(\vN\cdot w)=-\vN(\vN\cdot \vu).$$
 If $\vu$ belongs to a function space $E$ on which the Riesz transforms operate continuously, we find a particular solution $(\vv,\vw)$ by the formulas
 $$ \vv= \mathcal{R}\wedge(\mathcal{R}\wedge\vu) \text{ and } \vw=-\mathcal{R}(\mathcal{R}\cdot\vu) $$ where the vectorial Riesz transform is given as $$\mathcal{R}= \frac 1{\sqrt{-\Delta}} \vN.$$
 If $E$ contains no other harmonic function than the null function, then this decomposition $\vu=\vv+\vw$ is unique.
 
 The operator $\vu\mapsto  \mathcal{R}\wedge(\mathcal{R}\wedge\vu)$ is called the Leray projection operators (for $E=L^2$, it is the orthogonal projection of square-integrable vector fields on divergence-free square integrable vector fields) and is usually written as $\mathbb{P}$. This allows to get rid of the pressure in the Navier--Stokes equations and to rewrite the system as
 $$ \vu=\Delta\vu-\mathbb{P}((\vu\cdot\vN)\vu)$$ with $\vu(0,.)=\vu_0$ where $\vN\cdot\vu_0=0$.
 
 This way of eliminating the pressure $p$  (or expressing it as a function of the velocity $\vu$ by the formula $\vN p= \mathcal{R}(\mathcal{R}\cdot(\vu\cdot\vN)\vu)$) is quite general, and is applied to the study of (weak) solution $\vu$ in a large variety of function spaces. The justification for such computations has been given for instance by Furioli, Lemari\'e--Rieusset and Terraneo in the case of uniformly square-integrable solutions (vanishing at infinity)  \cite{FLT00, Lem02} or recently by Fernandez-Dalgo and Lemari\'e-Rieusset in the case of locally square integrable solutions with low increase at infinity \cite{FLR19}.

The nature of the Leray projection operator has a deep impact on the properties of solutions to the Navier--Stokes equations. Main features of the convolution kernel of the operator are that the kernel is not compactly supported, meaning that the operator is non-local and involves integration over the whole space, and that it has a slow decay at infinity (as $\mathbb{P}$ has a kernel homogeneous of degree $-3$, the kernel decays only as $\vert x\vert^{-3}$ and its derivatives as $\vert x\vert^{-4}$).  Writing, for a divergence-free vector field $\vu$,
$$   \mathcal{R}(\mathcal{R}\cdot(\vu\cdot\vN)\vu)= \vN( (\mathcal{R}\otimes\mathcal{R})\cdot (\vu\otimes\vu)),$$
Dobrokhotov and Shafarevich \cite{DSh94} proved that the spatial decay at infinity of ``rapidly'' decaying solutions was governed by the kernels $\partial_j\partial_k\partial_l G$  of the operators $ \partial_j R_kR_l$ (where $G$ is the Green function, fundamental solution of  the Laplacian operator : $G(x)=\frac 1{4\pi\vert x\vert}$, $(-\Delta) G=\delta$). More precisely, 
if $
 \lim_{x\rightarrow\infty}  \vert x\vert^4 \vert \vu^0(x)\vert=0 $ (as it is the case, for instance, for the Millenium Prize problem) and if  $(\vu,p)$ is a  classical solution of the Navier--Stokes problem on a strip $[0,T]\times\R^3$, then , for $0<t<T$,
 $$ u(t,x)= -\sum_{j=1}^3 \sum_{l=l}^3 d_{j,l}(t) \vN \partial_j\partial_l G(x)+o(\vert x\vert^{-4})
$$
with $d_{j,l}(t)=\int_0^t \int u_j(s,x) u_l(s,x)\, dx\, ds$. This means that the good decay of $\vu^0$  (as $o(\vert x\vert^{-4}$) is instantaneously lost whenever one of the integrals $\int u_j^0 u_l^0\, dx$ (with $j\neq l$) or $\int( u^0_j(x))^2 -(u^0_l(x))^2\, dx$ (with $j\neq l$) is not equal to $0$; in that case,  we have $\liminf_{x\rightarrow +\infty} \vert x\vert^4 \vert \vu(t,x)\vert >0$ for $t$ close enough to $0$. This instantaneous spreading has been studied by Brandolese and Meyer in \cite{BrM02}.
 
\subsection*{Lebesgue--Gevrey estimates.}
A less direct application of singular integrals to the study of the Navier--Stokes equations can be found in the treatment of Gevrey regularity of solutions in the Lebesgue space $L^3(\R^3)$ that has been proposed by Lemari\'e-Rieusset \cite{Lem00, Lem04}.

The idea starts from the result of Kato on existence of solutions for initial data in $L^3$ \cite{Kat84}.  One transforms the differential problem
$$ \partial_t \vu=\Delta \vu-\mathbb{P} \vN\cdot (\vu\otimes \vu),\quad \vu(0,.)=\vu^0
$$  into an integro-differential problem by solving
$$ \vu=e^{t\Delta} \vu^0-\int_0^t e^{(t-s)\Delta} \mathbb{P} \vN\cdot (\vu\otimes \vu)\, ds= e^{t\Delta} \vu^0-B(\vu,\vu)$$ where the bilinear operator $B$ is defined as
\begin{equation}\label{bilinear}B(\vu,\vv)(t,.)=\int_0^t e^{(t-s)\Delta} \mathbb{P} \vN\cdot (\vu\otimes \vv)\, ds. \end{equation}  By the contraction principle, if $B$ is bounded on a Banach space $\mathbb{E}_T$ (with operator norm $C_{\mathbb{E}_T}$) of functions defined on $(0,T)\times \R^3$, then, for $\vu^0$ small enough ($\|e^{t\Delta}\vu^0\|_B<\frac 1 {4C_{\mathbb{E}_T}}$), one can find a solution $\vu\in \mathbb{E}_T$. Now, if $\vu^0\in L^3$, we have 
\begin{equation}\label{injlebbesov} \sup_{t>0} t^{1/4} \|e^{t\Delta}\vu^0\|_6<+\infty\text{ and } \lim_{t\rightarrow 0}  t^{1/4} \|e^{t\Delta}\vu^0\|_6=0.\end{equation}
On the other hand, for every $t>0$,  the operator $e^{ t\Delta} \mathbb{P} \vN\cdot $ is given by convolutions with kernels $e^{t\Delta} \partial_j\partial_k\partial_l G$  that are in $L^1$ with 
$$\| e^{t\Delta} \partial_j\partial_k\partial_l G\|_1\leq C \frac 1{\sqrt t}. $$ We then use  the regularizing properties of the heat kernel in Lebesgue spaces : for $1\leq p\leq q$ and for $t>0$
$$ \|e^{t\Delta} f\|_q\leq C_{p,q}  t^{\frac 3 2 (\frac 1 q -\frac 1 p)} \|f\|_p.
$$ We then have
\begin{equation*}\begin{split} \|B(\vu,\vv)\|_6 \leq& \int_0^t \|  e^{\frac {t-s}2\Delta} \left( e^{\frac {t-s}2\Delta} \mathbb{P} \vN\cdot (\vu\otimes \vv)\right)\|_6\, ds \\ \leq &C \int_0^t \frac 1{(t-s)^{1/4}}  \frac 1{(t-s)^{1/2}}  \frac 1 {s^{1/2}} \|s^{1/4}\vu(s,.)\|_6 \|s^{1/4}\vv(s,.)\|_6\, ds. \end{split}\end{equation*}
It means that $B$ is bounded on
$$ \mathbb{E}_T=\{\vu(t,x)\ /\  \sup_{0<t<T} t^{1/4} \|e^{t\Delta}\vu^0\|_6<+\infty\text{ and } \lim_{t\rightarrow 0}  t^{1/4} \|e^{t\Delta}\vu^0\|_6=0\}$$ (with an operator norm that does not depend on $T$). Moreover, a solution in $\mathbb{E}_T$ will satisfy
$$   \|B(\vu,\vu)\|_3 \leq \int_0^t \|  e^{\frac {t-s}2\Delta} \left( e^{\frac {t-s}2\Delta} \mathbb{P} \vN\cdot (\vu\otimes \vu)\right)\|_3\, ds\leq C \int_0^t  \frac 1{(t-s)^{1/2}}  \frac 1 {s^{1/2}} (\|s^{1/4}\vu\|_6)^2\, ds$$ so that $\vu\in L^\infty ((0,T), L^3)$. (As a matter of fact, one even finds that $\vu\in \mathcal{C}([0,T), L^3)$).

Now, if we want to mimick the proof of the Gevrey regularity we saw for Fourier-Herz spaces, one must factor out in the Fourier--Navier--Stokes equations a term $e^{-\sqrt t \|\xi\|}$ and control the action of the factor $e^{\sqrt{t}\|\xi\|-\sqrt s \| \eta\|-\sqrt \|\xi-\eta\|}$. It is not enough to control the size of the factor, as we are dealing now with Fourier transforms of functions in $L^3$ or in $L^6$, that are no longer functions but singular distributions. The control is then given by the theory of singular integrals, and more precisely of Marcinkiewicz multipliers (as described in \cite{Ste71} for instance).  More precisely, we factor out in the Fourier transform a term of the form $e^{\sqrt t \|\xi\|_1}$, where $\|\xi\|_1=\vert\xi_1\vert  + \vert\xi_2\vert+\vert\xi_3\vert$. We shall write $e^{-\sqrt t D_1}$ for the convolution operator with symbol $e^{-\sqrt t \|\xi\|_1}$ and $e^{\sqrt t  D_1}$ for the convolution operator with symbol $e^{\sqrt t \|\xi\|_1}$. We then have to study the equation for $e^{\sqrt tD_1}\vu=\vU$ which is given by
\begin{equation*}\begin{split} \vU= e^{\frac t 2\Delta}&  \left(e^{\frac t 2\Delta} e^{\sqrt t D_1}\right) \vu_0\\&- \int_0^t e^{\frac{(t-s)}2\Delta}\mathbb{P} \vN\cdot \left( e^{\frac{(t-s)}2\Delta} e^{\sqrt{t-s} D_1}    e^{(\sqrt t-\sqrt{t-s} -\sqrt s)D_1}\big( e^{\sqrt s D_1} (e^{-\sqrt s D_1}\vU\otimes e^{-\sqrt sD_1}\vU)\big)\right) \, ds.\end{split}\end{equation*}
The operator $e^{\frac t 2 \Delta}e^{\sqrt t D_1}$ is a tensor product of one-dimensional convolution operators associated to Marcinkiewicz multipliers $e^{-\frac t 2\xi_j^2+\sqrt t \vert \xi_j\vert}$. Similarly, the operator  $e^{(\sqrt t-\sqrt{t-s} -\sqrt s)D_1} $  is a tensor product of one-dimensional convolution operators associated to Marcinkiewicz multipliers
 $e^{ (\sqrt t-\sqrt {t-s}-\sqrt s)\vert \xi_j\vert}$. The bilinear operator
$$ T(f,g)=e^{ \sqrt sD_1} (e^{-\sqrt s D_1}f \times e^{-\sqrt sD_1} g)$$ can similarly be written as a sum of  tensor products of one-dimensional convolution operators associated to Marcinkiewicz multipliers : if $S_j$ is associated to the multiplier $1_{\xi_j>0}$, $T_j$ to the multiplier $1_{\xi_j}<0$ and $Z_j$ to the multiplier $e^{- \sqrt s \vert \xi_j\vert}$ and if $W_j$ is the unbounded operator associated to the multiplier  $e^{+\sqrt s \vert\xi_j\vert}$, then, for $f_j$, $g_j\in L^p(\R)$,
\begin{equation*} \begin{split}W_j(Z_jf_j\times Z_jg_j)=&S_jf_j\times S_jg_j+ T_jf_j\times T_jg_j+S_j(S_jf\times Z_j^2T_jg)\\ &+S_j(Z_j^2T_jf\times S_jg) + T_j(Z_j^2S_jf\times T_jg)+T_j(T_jf\times Z_j^2S_jg).
\end{split}\end{equation*}
Thus, using the contraction principle, we find that, if $\vu_0\in L^3$, we have a solution $\vu=e^{-\sqrt t D_1} \vU$ of the Navier--Stokes equations  on a small enough time interval $(0,T)$ such that $\sup_{0<t<T} t^{1/4} \|\vU\|_6<+\infty$ and $\vU\in L^\infty((0,T),L^3)$.

\subsection*{Maximal regularity for the heat kernel.}
Another way of using singular integrals for the study of solutions to the Navier--Stokes equations is the proof proposed  by Monniaux in \cite{Mon99} for the uniqueness of solutions in  $\mathcal{C}([0,T), L^3)$. (Local) existence of solutions in $L^3$ (for an initial value $\vu^0\in L^3$) had been proved by Kato in 1984 \cite{Kat84}, but uniqueness remained open until 1997, when Furioli, Lemari\'e-Rieusset and Terraneo \cite{FLT00} proved uniqueness by using Besov spaces.

The proof by Monniaux is very simple. If $\vu$ is a solution in $\mathcal{C} ([0,T),L^3)$ and $\vu_K$ is the solution provided by Kato in $\mathcal{C (}[0,T),L^3)$ with the additional property that $\lim_{t\rightarrow 0} t^{1/4} \|\vu_K\|_6=0$, then the function $\vw=\vu-\vu_K$ satisfies the identity
$$ \vw=-B(\vu_K,\vw)-B(\vw,\vu_K)-B(\vw,\vw).$$
 We shalll write that $\vw$ is an eigenvector of the linear transform
 $$ \vv \mapsto L(\vv) =-B(\vu_K,\vv)-B(\vv,\vu_K)-B(\vw,\vv).$$
 We want to estimate $L(\vv)$ in $L^3((0,S),L^3)$, for $S<T$. We have
 $$ \|B(\vu_K,\vv)(t,.)\|_3\leq C \int_0^t  \frac 1{(t-s)^{1/4}}   \frac 1{\sqrt{t-s}} \frac 1 {s^{1/4}} \|\vv(s,.)\|_3  \|s^{1/4}\vu_K(s,.)\|_6 \, ds$$ so that, since multiplication is bounded from $L^{4,\infty}\times L^3$ to $L^{12/7,3}$ and convolution is bounded from $L^{12/7,3}\times L^{4/3,\infty}$ to $L^{3,3}=L^3$, 
$$  \|B(\vu_K,\vv)\|_{L^3((0,S),L^3)} \leq C \|\vv\|_{L^3((0,S),L^3)} \sup_{0<s<S} s^{1/4} \|\vu_K(s,.)\|_6.$$ Similarly, we have
$$  \|B(\vv,\vu_K)\|_{L^3((0,S),L^3)} \leq C \|\vv\|_{L^3((0,S),L^3)} \sup_{0<s<S} s^{1/4} \|\vu_K(s,.)\|_6.$$ 
 For estimating $B(\vw,\vv)$, we write
 $$ \partial_j R_kR_l(w_kv_l)=-\Delta  R_jR_kR_l \frac 1{\sqrt{-\Delta}} (w_kv_l$$
 and use the inequality on Riesz potential
 $$ \|\frac 1{\sqrt{-\Delta}}(f)\|_3\leq \|f\|_{3/2}.$$ Thus, we find that
 $$ \frac 1 \Delta (\mathbb{P}\vN\cdot(\vw\otimes\vv))\in  L^3 L^3.
$$
Maximal regularity in $L^3 L^3$ for the heat kernel states that
$$ \|\int_0^t e^{(t-s)\Delta} \Delta f\, ds\|_{L^3((0,S), L^3)}\leq C \|f\|_{L^3((0,S),L^3)}$$ where the constant $C$ does not depend on $S$. Thus, we have
$$ \|B(\vw,\vv)\|_{L^3((0,S),L^3)} \leq C \|\vv\|_{L^3((0,S),L^3)} \|\vw\|_{L^\infty((0,S),L^3)} $$ By continuity of $\vw$ in $L^3$, we find that $\lim_{S\rightarrow +\infty}  \|\vw\|_{L^\infty((0,S),L^3})=0$. Thus, for $S$ small enough, ${L}$ is contractive on $L^3((0,L^3),L^3)$. Hence, the fixed point $\vw$ is equal to $0$, and $\vu=\vu_K$ on $(0,S)$. The end of the proof follows by a bootstrap argument.

The maximal regularity property is linked to singular integrals, but no longer on $\R^3$ but on the space $ \R\times\R^3$ endowed with the parabolic distance $\delta((t,x),(s,y))=\sqrt{\vert t-s\vert +\vert x-y\vert^2}$. Together, with the Lebesgue measure on $X=\R\times\R^3$, $\delta$ provides $X$ with a structure of homogeneous space (as studied by Coifman and Weiss) \cite{CoW71}. We have, for $f$ supported in $[0,+\infty)\times\R^3$, 
$$1_{t>0} \int_0^t e^{(t-s)\Delta} \Delta f\, ds=\iint_X K(t-s,x-y) f(s,y)\, ds\, dy$$ where $K$ is a convolution operator associated to the Fourier multiplier 
$$ m(\tau,\xi)=-\frac{\xi^2}{\xi^2+i\tau}$$ (where we consider the Fourier transform $\mathcal{F}_{t,x}f(\tau,\xi)=\iint_X f(s,y) e^{-i(t\tau+x\cdot\xi)}\, dt\, dx$). We have
$$ \sup_{\alpha\in \N_0^3, \beta\in \N_0} \sup_{(\tau,\xi)\neq (0,0)}   (\vert \tau\vert^{1/2}+\vert\xi\vert)^{\vert \alpha\vert+2\beta} \vert \frac{\partial^\alpha}{\partial\xi^\alpha} \frac{\partial^\beta}{\partial\tau^\beta}m(\tau,\xi)\vert <+\infty.$$
Thus, $m$ can be seen as a Marcinkiewicz multiplier on the parabolic space $\R\times\R^3$.
`
\subsection*{Marcinkiewicz multipliers for bilinear operators.}
In 1978, Coifman and Meyer extended the theory of multipliers to the setting of bilinear operators \cite{CoM78, CoM91}.  They consider a smooth function $\sigma$ on $\R^d\times \R^d$ such that, for all  $\alpha$, $\beta\in \N_0^d$,
$$ \sup_{(\xi,\eta)\neq (0,0)} (\vert \xi\vert+\vert \eta\vert)^{\vert\alpha\vert+\vert\beta\vert} \left\vert \frac{\partial^\alpha}{\partial\xi^\alpha}\frac{\partial^\beta}{\partial \eta^\beta}\sigma(\tau,\xi)\right\vert<+\infty$$
and they define
$$ T_\sigma(f,g)=\frac 1{(2\pi)^{2d}}\iint e^{i(\xi+\eta)\cdot x} \mathcal{F}_xf(\xi) \mathcal{F}_xg(\eta)\, d\xi\, d\eta.$$
$T$ is bounded from $L^\infty\times L^p$ to $L^p$ for every $1<p<+\infty$. One key property is that, for fixed $f\in L^\infty$, $g\mapsto T(f,g)$ is not a convolution operator but  is a generalized Calder\'on--Zygmund operator (in the sense of \cite{CoM78b}).

This theory has been applied by Kato and Ponce to derive a useful commutator estimate that they applied to the study of the regularity of solutions to the Navier--Stokes equations or to the Euler equations \cite{KaP88}.

\section{The Hardy--Littlewood maximal function.}

\subsection*{Kato's mild solutions and maximal functions.}

 Our treatment of the Navier--Stokes equations through the cheap Navier--Stokes equation was very elementary, using absolute values and convolution inequalities in the frequency variables. C. Calder\'on  \cite{Cal93} noticed that we can deal with  the equations in the space variable in an equivalently elementary way, through the use of the maximal function, another basic tool in harmonic analysis introduced by Hardy and Littlewood in 1930 \cite{HaL30}. We shall write $\mathcal{M}_f$ for the maximal function of $f$ :
 $$ \mathcal{M}_f(x)=\sup_{r>0} \frac 1{\vert B(x,r) \vert} \int_{B(x,r)} \vert f(y\vert\, dy.$$

 A basic result on maximal functions \cite{Gra08} is the control of convolution with radial kernels. More precisely, if a function $g$ admits a majorant $k$ ($\vert g(x)\vert\leq k(x)$) such that $k$ is integrable, radial and radially non-increasing, then 
 $$ \vert g*f(x)\vert \leq \|k\|_1 \mathcal{M}_f(x).$$
 
 In order to deal with
 the integro-differential problem 
$$ \vu=e^{t\Delta} \vu^0-\int_0^t e^{(t-s)\Delta} \mathbb{P} \vN\cdot (\vu\otimes \vu)\, ds= e^{t\Delta} \vu^0-B(\vu,\vu),$$
Calder\'on writes
$$ e^{t\Delta} \vu^0=e^{t\Delta}  (-\Delta)^{1/4} (-\Delta)^{-1/4}\vu^0$$
and
$$ B(\vu,\vv)=\int_0^t  e^{\frac{(t-s)}2 \Delta} \mathbb{P} \vN\cdot  e^{\frac{(t-s)}2\Delta}  (-\Delta)^{1/4} (-\Delta)^{-1/4}(\vu\otimes \vv)\, ds.$$ and then uses the inequalities
\begin{itemize}
\item[$\bullet$] $\vert e^{t\Delta}f (x)\vert\leq C \int \frac {\sqrt t}  {(\sqrt t+\vert x-y\vert)^4} \vert f(y)\vert\, dy$
\item[$\bullet$] $\vert e^{t\Delta} (-\Delta)^{1/4} f (x)\vert\leq C \int \frac 1 { (\sqrt t+\vert x-y\vert)^{7/2}} \vert f(y)\vert\, dy$
\item[$\bullet$] $\vert e^{t\Delta} (-\Delta)^{1/4} R_jR_k\partial_l f (x)\vert\leq C \int \frac 1 { (\sqrt t+\vert x-y\vert)^{4}} \vert f(y)\vert\, dy$
\item[$\bullet$] $\vert (-\Delta)^{-1/4} (fg)(x)\vert \leq C \int \frac 1 {\vert x-y\vert^{5/2}} \vert f(y) g(y)\vert\, dy$
\end{itemize}
The last inequality is just a consequence of the correspondence of fractional integration $(-\Delta)^{-\alpha/2}$ ($0<\alpha<3$) with the Riesz potentials $I_\alpha$ :
$$(-\Delta)^{-\alpha/2} f(x)= I_\alpha(f)(x)=c_\alpha\int_{\R^3} \frac 1{\vert x-y\vert^{3-\alpha}}\, f(y)\, dy.$$

With the first two inequalities, we get that 
$$ \sup_{t>0} \vert e^{t\Delta} \vu^0(x)\vert \leq  C \mathcal{M}_{\vert\vu_0\vert}(x)$$ and
$$\sup_{t>0} t^{1/4} \vert e^{t\Delta} \vu^0(x)\vert \leq  C \mathcal{M}_{I_{1/2}(\vert\vu_0\vert)}(x).$$ 
In particular, if $\vu^0$ belongs to $L^3$ then (as $I_{1/2}$ maps $L^3$ to $L^6$) we have $e^{t\Delta}\vu^0\in \mathbb{E} $ where
$$ f\in \mathbb{E}\Leftrightarrow \sup_{t>0} f(t,x)\in L^\infty(\R^3)\text{ and } \sup_{t>0} t^{1/4} f(t,x)\in L^6(\R^3).$$ Moreover, as $L^6\cap L^3$ is dense in $L^3$, we have
$$ \lim_{T\rightarrow 0} \| \sup_{0<t<T} t^{1/4}  e^{t\Delta}\vu^0\|_6=0.$$
Now if $\vu\in \mathbb{E}_T$ and $\vv\in \mathbb{E}_T$ where
$$ f\in \mathbb{E}_T \Leftrightarrow A_T(f)=\sup_{0<t<T}  t^{1/4} \vert f(t,x)\vert\in L^6,$$ we find (using the inequalities on $e^{t\Delta}(-\Delta)^{1/4}$ and on $e^{t\Delta}R_jR_k\partial_l$, the inequalities (for $0<t<T$)
$$ \vert B(\vu,\vv)(t,x)\vert \leq  C\int_0^t \frac 1{(t-s)^{3/4}} \frac 1 {\sqrt s} \mathcal{M}_{I_{1/2}(A_T(\vu) A_T(\vv))}(x)\, ds \leq C' t^{-\frac 1 4} \mathcal{M}_{I_{1/2}(A_T(\vu) A_T(\vv))}(x)$$
and 
$$ \vert B(\vu,\vv)(t,x)\vert \leq  C\int_0^t \frac 1{(t-s)^{1/2}} \frac 1 {\sqrt s} \mathcal{M}_{A_T(\vu) A_T(\vv)}(x)\, ds \leq C' \mathcal{M}_{A_T(\vu) A_T(\vv)}(x).$$
The first inequality gives that $B(\vu,\vv)$ still belongs to $\mathbb{E}_T$, and thus, if $T$ is small enough (to grant that $e^{t\Delta}\vu^0$ is small in $\mathbb{E}_T$), we find a solution $\vu$ to the Navier--Stokes problem; the second inequality gives us a control in $L^3$ nborm for this solution $\vu$. Thus, we recover a Kato-type solution
$\vu$ such that
$$ \sup_{0<t<T} \vert\vu(t,.)\vert\in L^3\text{ and } \sup_{0<t<T} t^{1/4} \vert\vu(t,.)\vert \in L^6.$$
The main difference with Kato's formalism is that, now,  we first take the supremum on $t$ before integrating in $x$.\\

 When $\vu^0$ is small in $L^3$, we have an even simpler proof of existence of a mild solution $\vu$ such that $\sup_{t>0}  \vert\vu(t,.)\vert$ belongs to $L^3$.  This result of Calder\'on is based on the fact that the bilinear operator $B$, which is not bounded on $L^\infty_t L^3_x$ \cite{Oru98}, is actually bounded on $L^3_x L^\infty_t$.                                                   If $\vu\in L^3_x L^\infty_t$ and $\vv \in L^3_x L^\infty_t$, then
 \begin{equation*}\begin{split}\vert B(\vu,\vv)(t,x)\vert \leq& C \int_0^t\int_{\R^3} \frac 1 { (\sqrt t+\vert x-y\vert)^{4}} \vert \vu(s,y)\vert\,  \vert\vv(s,y)\vert\, dy\\ \leq& C \int_{\R^3} \sup_{s>0} \vert \vu(s,y)\vert \sup_{s>0}\vert\vv(s,y)\vert \left[ \int_0^t  \frac 1 { (\sqrt t+\vert x-y\vert)^{4}} \, ds\right]\, dy\\ =& C'  \int_{\R^3}  \frac 1{\vert x-y\vert^2}\sup_{s>0} \vert \vu(s,y)\vert \sup_{s>0}\vert\vv(s,y)\vert \, dy
 \\ =& C'' I_{1/2}(  \sup_{s>0}\ \vert \vu(s,y)\vert \sup_{s>0}\vert\vv(s,y)\vert )(x).
 \end{split}\end{equation*}
 Thus, if $\mathcal{M}_{\vert \vu^0\vert} \leq U^0$ and if $U$ is a solution of the cheap equation 
 $$ U=U^0+ C''  I_{1/2}(U^2)$$ (solution which exists if $U^0$ is small enough in $L^3$), then we have a  solution $\vu$ with $\vert \vu(t,x) \vert \leq U(x)$.                                            

\subsection*{Hardy spaces and molecules.}

 We can rewrite the Hardy--Littlewood maximal function as 
 $$ \mathcal{M}_f(x)=\sup_{t>0} K_t*\vert f\vert(x)$$ where $$ K(x)=\frac 1{\vert B(0,1)\vert} \1_{B(0,1)}(x) \text{ and } K_t(x)=\frac 1 {t^3} K(\frac x t).$$ One can see clearly the role of scaling in this definition of the operator. Basic features for this operator are the boundedness on $L^p$ for $1<p\leq +\infty$
 $$ \|\mathcal{M}_f\|_p \approx \|f\|_p$$ and the lack of control in $L^1$ norm :
 $$ f\neq 0\implies \|\mathcal{M}_f\|_1=+\infty.$$
 
 The theory of Hardy spaces developed by Fefferman and Stein \cite{FeS72, Ste93} involves a modified maximal function : taking $\Phi\in\mathcal{S}$ a radially non-increasing smooth function, and defining $\Phi_t(x)=\frac 1 {t^3} \Phi(\frac x t)$, one defines
 $$ \mathcal{M}_f^{[\Phi]}(x)=\sup_{t>0} \vert \Phi_t*f(x)\vert.$$ From the properties that $\mathcal{M}_f^{[\Phi]}(x)\leq \mathcal{M}_f(x)$ and that $\lim_{t\rightarrow 0} \Phi_t*f=f$ in $\mathcal{S}'$ for every distribution $f\in\mathcal{S}'(\R^3)$, one finds that we have again, for $1<p\leq +\infty$, 
  $$ \|\mathcal{M}^{[\Phi]}_f\|_p \approx \|f\|_p.$$ 
  But, now, it turns out that there are many distributions $f$ such that $\mathcal{M}^{[\Phi]}_f$ is integrable or belongs to $L^p$ for some $p\in (0,1)$. The Hardy space $\mathcal{H}^p$ is defined for $0<p<+\infty$ by the property :
  $$ f\in\mathcal{H}^p\Leftrightarrow \mathcal{M}^{[\Phi]}_f\in L^p.$$
  
  An important feature of Hardy spaces is their duality property with ${\rm BMO}$ or with homogeneous H\"older spaces : the dual of $\mathcal{H}^1$ can be identified with ${\rm BMO}$ and the dual of $\mathcal{H}^p(\R^3)$ for $\frac 3 4<p<1$ can be identified with the homogeneous H\"older space $\dot B^\alpha_{\infty,\infty}$ with $\alpha=\frac 3 p-1$. This has been used by Kozono and Taniuchi  \cite{KzT00} to prove weak-strong uniqueness solutions when  the Navier--Stokes problem with initial value $\vu^0\in L^2$ generates a weak Leray solution in $L^\infty_t L^2_x\cap L^2_t \dot H^1_x$ (with Leray energy inequality) and a solution in  $L^\infty_t L^2_x\cap L^2_t \dot H^1_x\cap L^2_t {\rm BMO}_x$ : the proof relies on the proof by Coifman, Lions, Meyer and Semmes \cite{CLMS92} that, for a vector field $\vu\in L^2$  that is divergence-free ($\vN\cdot\vu$) and a vector field $\vv\in L^2$ that is curl-free ($\vN\wedge\vv=0$), we have $\vu\cdot\vv\in\mathcal{H}^1$. Thus, the usual estimate for weak-strong uniqueness
  $$ \|\vv-\vu\|_2^2+2\int_0^t \|\vN\otimes (\vu-\vv)\|_2^2\, ds \leq 2 \int_0^t\int_{\R^3} \vu\cdot\left( (\vu-\vv)\cdot \vN(\vu-\vv)\right)\, dx\, ds$$
 is turned to
   $$ \|\vv-\vu\|_2^2+2\int_0^t \|\vN\otimes (\vv-\vu)\|_2^2\, ds \leq  C \int_0^t  \|\vv-\vu\|_2 \|\vN(\vv-\vu)\|_2 \|\vu\|_{\rm BMO}\, \, ds$$ which leads to a Gronwall estimate.
   
   Another important feature of Hardy spaces is their atomic decomposition, as described for instance in \cite{CoW71}.  For $3/4<p\leq 1$, we have that a distribution $f$ belongs to $\mathcal{H}^p(\R^3)$ if and only if it can be written as a sum $f=\sum_{j\in\N} \lambda_j  a_j$ where $\sum_{j\in\N} \vert \lambda_j\vert^p<+\infty$ and $a_j$ is a $\mathcal{H}^p$ atom : there exists some $r_j>0$ and some 
 $x_j\in \R^3$ such that $a_j$ is supported in the ball $B(x_j,r_j)$, $\|a_j\|_2\leq \vert B(x_j,r_j)\vert^{-\frac 2 {2-p}}$ and $\int a_j\, dx=0$.
 
 Atoms are not stable under the action of Calder\'on--Zygmund convolution operators with a non-local singular kernel, because compactness of supports is destroyed by the convolution. But if we relax the conditions on $a_j$ into, for some $r_j>0$ and $x_j\in\R^3$,   $\|a_j\|_2\leq \vert B(x_j,r_j)\vert^{-\frac 2 {2-p}}$, $\| \ \vert x-x_j\vert\,  a_j\|_2\leq r_j  \vert B(x_j,r_j)\vert^{-\frac 2 {2-p}}$ and $\int a_j\, dx=0$ [$a_j$ is no longer an atom for $\mathcal{H}^p$, but it is called a molecule], the situation is much better. If $T$ is a convolution operator with a  Calder\'on--Zygmund kernel, then there exists a constant  $C>0$ such that the image $\frac 1 C T(a_j)$ of a molecule  is still a molecule (associated to the same center $x_j$ and the same radius $r_j$). 
 
 There are very few examples of the use of Hardy molecular decompositions in fluid mechanics. We may quote  a paper of Chamorro on advection-diffusion in the setting of a non-local diffusion and a rough drift \cite{Cha14}. Futioli and Terraneo \cite{FuT02} studied the Cauchy problem for the Navier--Stokes equations when the Laplacian of the initial value $\vu_0$ is a $\mathcal{H}^1$ molecule.
Their results were extended by Brandolese in \cite{Bra01}.

\subsection*{Wavelets.}
Atomic or molecular decompositions lead quite naturally to wavelets. However, the basic atoms that generate wavelet decompositions are usually more regular than simply Lebesgue measurable and are assumed to have some H\"oder regularity. In that case, one works more in the setting of Besov spaces than of Hardy spaces. A systematic approach of Besov space through atomic decompositions has been proposed by many authors, including the seminal paper of Frazier and Jawerth \cite{FrJ90}.

However, the first approach of the Navier--Stokes equations with a decomposition on wavelet bases was performed by Federbush   \cite{Fed93} in yet anothe space, the Morrey space $\dot M^{2,3}$ (see the subsection  on Morrey spaces in section \ref{FS}). The study by Federbush was based on the use of divergence-free vector wavelet bases  \cite{BaF95, Lem92}. Divergence-free vavelet bases were also used by Urban \cite{Urb95} for the numerical approximation of the equations of fluid mechanics.

There have been many claims that wavelet analysis of turbulent signals may provide valuable insights in the actual structure of turbulent fluids  \cite{FKPS99, FrZ93}, especially in the frame of self--similar universality laws  such as studied by Frisch \cite{Fri95}. But Meyer proved that the claim that wavelets were asymptotically decorrelated in the non-linearity of the Navier--Stokes equations was unfounded \cite{Mey99}.

\section{Function spaces} \label{FS}
Many function spaces of measurable or differentiable functions have close relationships with harmonic analysis, and their theory were developed quite extensively in the books of Stein :  Lorentz spaces in \textit{Introduction to Fourier Analysis on Euclidean
Spaces} \cite{StW71}, Besov spaces in  \textit{Introduction to Fourier Analysis on Euclidean
Spaces} \cite{Ste71}, $\rm BMO$, tent spaces  or Muckenhoupt weights in \textit{Harmonic Analysis} \cite{Ste93}. As a matter of fact, all those spaces are met in the modern study of Navier--Stokes equations developed in the 90's. More recently, use of Morrey spaces has been developed as well  by many authors (see \cite{Lem16} for references).

\subsection*{Lorentz spaces.}

Sobolev spaces $W^{k,p}$ of functions in $L^p$ such that their derivatives (in the sense of distributions) up to order $k$ are still in $L^p$ can ce extended for $1<p<+\infty$  to the scale of spaces $H^s_p$ defined, for $s\in\R$, by $$f\in H^s_p\Leftrightarrow f\in\mathcal{S}'\text{ and } \mathcal{F}_x^{-1} \left( (1+\vert\xi\vert^2)^{s/2} \mathcal{F}_x f \right)\in L^p.$$
For $1<p<+\infty$ and $k\in \N_0$, we have $W^{k,p}=H^k_p$ \cite{Ste71}. The Sobolev embeddings then state that, for $0\leq s< \frac 3 p$ (and $1<p<+\infty$) we have
$$ H^s_p \subset L^r\text{ with } \frac 1 r=\frac 1 p-\frac s 3.$$ The \textit{sharp Sobolev embedding} 
states more precisely that
$$ H^s_p \subset L^{r,p}\subset L^r \text{ with } \frac 1 r=\frac 1 p-\frac s 3$$ where $L^{r,p}$ is a Lorentz space. This can be done  through various methods that have been developed by Stein; for instance :
\begin{itemize}
\item[$\bullet$]  Let $\mathcal{J}^s$ be the convolution with the Bessel kernel associated with the Fourier multiplier $(1+\vert\xi\vert^2)^{-s/2}$; convolution with $\mathcal{J}^s$ maps $L^p$ onto $H^s_p$ for $1<p<+\infty$; r-the Sobolev embeddings state that it maps $L^p$ to  $L^r$ with $ r=\frac{3p}{3-sp}$ when $0\leq s<\frac 3 p$; then, picking $p_0$ and $p_1$ with $1<p_0<p<p_1<\frac 3 s $,  the Marcinkiewicz interpolation theorem (as extended by Stein and Weiss \cite{Mar39, StW71}) gives the boundedness of  $\mathcal{J}^s$ from $L^p$ to $L^{\frac {3p}{3-sp},p}$ as an interpolation of the boundedness of $\mathcal{J}^s$ from $L^{p_0}$ to $L^{\frac{3p_0}{3-sp_0}}$ and from $L^{p_1}$ to $L^{\frac{3p_1}{3-sp_1}}$.
\item[$\bullet$] For $0<s<3$, the kernel $K_s$ of the convolution operator $\mathcal{J}^s$ satisfies
$$ \vert K_s(x)\vert \leq C \frac 1{\vert x\vert^{3-s}}$$ and thus $K_s$ belongs to the Lorentz space $L^{\frac 3{3-s},\infty}$. Convolution in Lorentz spaces has been studied by  O'Neil \cite{ONe63} following ideas of Stein. In particular, we have $L^{p,q}* L^{r,s}\subset L^{t,u}$ with $\frac 1 t=\frac 1 p+\frac 1 r-1$ and $\frac 1 u =\min(\frac 1 q+\frac 1 s,1)$ (whenever $1< t< +\infty$). Applying this to $L^p=L^{p,p}$ and $L^{\frac 3{3-s},\infty}$ gives the desired embedding.
\end{itemize}

Due to their good properties of interpolation and to their simple convolution and product laws, Lorentz spaces have turned out to be very efficient tools for providing sharp estimates in Lebesgue norms. For instance, the Hardy inequality
  $$ \int_{\R^3}  \frac{ \vert f\vert^p} {\vert x\vert^{sp}}\, dx \leq C_{s,p} \int_{\R^3} \vert (-\Delta)^{s/2} f\vert^p\, dx$$ for $f\in H^s_p$, $1<p<+\infty$  and $0<s<3/p$ is a  direct consequence of the facts that the kernel of $(-\Delta)^{-s/2}$ (i.e. the Riesz potential $\mathcal{I}^s$) belongs to $L^{\frac 3{3-s},\infty}$, the multiplier $\frac 1{\vert x\vert^s}$ belongs to $L^{\frac 3 s,\infty}$, the convolution maps $L^p\times L^{\frac 3{3-s},\infty}$ to $L^{\frac {3p}{3-sp},p}$ and the pointwise product maps $L^{\frac {3p}{3-sp},p}\times L^{\frac 3 s,\infty}$ to $L^{p,p}=L^p$.
  
  Besides being a useful tool for refining inequalities, Lorentz spaces occur as a natural setting in various problems in the study of the Navier--Stokes equations, especially in problems with critical scaling. For instance, the bilinear operator $B$  defined by equation (\ref{bilinear}) is bounded on $L^\infty((0,T),L^p)$ for $p>3$, with a norm $$ \| B\|_{\mathcal{B}(L^\infty L^p\times L^\infty L^p\mapsto L^\infty L^p)}= C_p T^{\frac 1 2 (1-\frac 3 p)}.$$ But it is no longer bounded on $L^\infty L^3$  \cite{Oru98}. It turns out that, however, it is bounded on $L^\infty L^{3,\infty}$, as proved by Meyer  \cite{Mey99}.
  
  Kozono and Nakao \cite{KzN96} studied time-periodic solutions for the Navier--Stokes equations with a time-periodic forcing. They found solutions in $L^\infty L^{3,\infty}$. More precisely, one studies the equations 
$$ \partial_t \vu=\Delta \vu-\mathbb{P} \vN\cdot (\vu\otimes \vu)+\mathbb{P}\vN\cdot \mathbb{F}
$$  where the forcing tensor $\mathbb{F}$ is time-periodic, and one seeks for a solution $\vu$ which is still time-periodic. If we assume that $\mathbb{F}$ belongs to $L^\infty L^{3/2,\infty}$, then we define $\vU^0$ as $$\vU^0=\int_{-\infty}^t   e^{(t-s)\Delta} \mathbb{P} \vN\cdot \mathbb{F}\, ds$$
and we find that $\vU^0$ belongs to $L^\infty L^{3,\infty}$. Thus, looking for time-periodic solutions of the Navier--Stokes solutions with time-periodic tensor $\mathbb{F}$  is turned into the solving of 
the  integro-differential problem 
$$ \vu=\vU^0-\int_{-\infty}^t e^{(t-s)\Delta} \mathbb{P} \vN\cdot (\vu\otimes \vu)\, ds= \vU^0-B_\infty(\vu,\vu)$$ where the bilinear operator   $B_\infty$  is defined as
 $$ B_\infty(\vu,\vv)(t,.)=\int_{-\infty}^t e^{(t-s)\Delta} \mathbb{P} \vN\cdot (\vu\otimes \vv)\, ds.  $$
  As $B_\infty$ is bounded on $L^\infty L^{3,\infty}$, the Banach contraction principle will give us a solution as soon as $\mathbb{F}$ is small enough.

    Meyer \cite{Mey99}  applied the boundedness of $B$  on $L^\infty L^{3,\infty}$ to another problem, namely  uniqueness of solutions  of the Navier--Stokes equations in  $\mathcal{C}([0,T), L^3)$.  We already discussed this problem. 
   Let $\vu^0\in L^3$.  If $\vu$ is a solution in $\mathcal{C}[0,T),L^3)$ and $\vu_K$ is the solution provided by Kato in $\mathcal{C}[0,T),L^3)$ with the additional property that $\lim_{t\rightarrow 0} t^{1/4} \|\vu_K\|_6=0$, then the function $\vw=\vu-\vu_K$ satisfies the identity
$$ \vw=-B(\vu_K,\vw)-B(\vw,\vu_K)-B(\vw,\vw).$$
The main idea of the proof of uniqueness initially given by  Furioli, Lemari\'e-Rieusset and Terraneo \cite{FLT00} was to establish a contractive estimate (locally in time) on $\vw$ to prove that $\vw$ is equal to $0$. The difficult term is $B(\vw,\vw)$, as $B$ is not bounded on $L^\infty L^3$. But, as $B$ is bounded on $L^\infty L^{3,\infty}$, it is easy to prove that, for $0<S<T$, 
$$ \sup_{0<t<S} \|\vw(t,.)\|_{L^{3,\infty}}\leq   C  \sup_{0<t<S} \|\vw(t,.)\|_{L^{3,\infty}} \left( \sup_{0<s<S} s^{1/4} \|\vu_K(s,.)\|_6+  \sup_{0<s<S} \|\vw(t,.)\|_{L^{3,\infty}} \right).$$ 
By continuity of both $\vu$ and $\vu_K$ in $L^3$ norm, and by the embedding $L^3\subset L^{3,\infty}$, we have  $\lim_{t\rightarrow 0} \|\vw(t,.)\|_{L^{3,\infty}}=0$ and we find that we have a contractive estimate for $\vw$ if $S$ is small enough.  Hence, the fixed point $\vw$ is equal to $0$, and $\vu=\vu_K$ on $(0,S)$. The end of the proof follows by a bootstrap argument.

    Another example where one naturally deals with Lorentz spaces is the study of self-similar solutions.
If $\vu$ and $p$ are solutions of 
$$ \partial_t\vu+(\vu\cdot\vN)\vu=\Delta\vu -\vN p,$$
$$  \vN\cdot \vu=0,$$
$$\vu(0,.)=\vu^0,$$
then, for $\lambda>0$, defining $\vu_\lambda(t,x)=\lambda \vu(\lambda^2 t,\lambda x)$, $p_\lambda(t,x)=\lambda^2 p(\lambda^2 t,\lambda x)$ and $\vu^0_\lambda(x)=\lambda\vu^0(\lambda x)$, we find that $\vu_\lambda$ and $p_\lambda$ are solutions of 
$$ \partial_t\vu_\lambda+(\vu_\lambda\cdot\vN)\vu_\lambda=\Delta\vu_\lambda -\vN p_\lambda,$$
$$  \vN\cdot \vu_\lambda=0,$$
$$\vu_\lambda(0,.)=\vu_\lambda^0.$$
Thus, provided that $\vu_0$ is homogeneous (so that $\vu^0_\lambda=\vu^0$), one may look for self-similar solutions (such that $\vu_\lambda=\vu$ and $p_\lambda=p$). However, if $\vu^0$ is homogeneous and is not equal to $0$, it cannot belong to the usual spaces (Lebesgue spaces $L^p$ or Sobolev spaces $H^s$), by lack of integrability either at $x=0$ or at $x=\infty$. B ut $L^{3,\infty}$ contains non-trivial homogenous functions, so that  the problem of looking for self-similar solutions is meaningful in the setting of this Lorentz space (this has been done by Barraza in 1996
  \cite{Bar96}).

\subsection*{Besov spaces.}
Besov spaces are usually viewed as the main tool of real variable harmonic analysis methods for the Navier--Stokes equations \cite{Lem02, Can04, BCD11}. However, the role played by Besov spaces has various aspects.

The most obvious occurence of Besov spaces, more related to the classical theory of parabolic equations than to real methods in harmonic analysis, is linked to the analysis of the heat kernel and the thermic characterization of Besov spaces. If $s<0$, then a distribution $f\in\mathcal{S}'$ will belong to $B^s_{p,q}$ ($1\leq p,q\leq +\infty$) if and only if $t^{s/2} \| e^{t\Delta}f\|_q$ belongs to $L^p((0,T), \frac {dt}t)$ (where $0<T<+\infty$). If $T=+\infty$, then $f$ belongs to the (realization of)  homogeneous  Besov space $\dot B^s_{p,q}$. Thus, the result of Fabes, Jones and Rivi\`ere \cite{FJR72} that $B$ is bounded on $L^p((0,T), L^q)$ when $T$ is finite, $3<q<+\infty$  and $\frac 2 p+\frac 3 q\leq 1$, or when $T=+\infty$, $3<q<+\infty$  and $\frac 2 p+\frac 3 q= 1$, implies that one may find a (local in time) solution $\vu\in L^p L^q$ to the Cauchy problem for the Navier-Stokes equations with initial value $\vu_0$ if and only $\vu_0\in B^{-\frac 2 p}_{q,p}$; this solution will be global if $\frac 2 p+\frac 3 q=1$ and $\vu_0$ is small enough in $\dot B^{-\frac 2 p}_{q,p}$.

Similarly, the inequality (\ref{injlebbesov}) we used to construct Kato solutions to the problem with $\vu_0\in L^3$ can be seen as the embedding $L^3\subset \dot B^{-\frac 1 2}_{6,\infty}$.

Homogeneous Besov spaces occur naturally in the study of the Navier--Stokes equations, due to the scaling invariance of the equations. If we want to study the initial value problem in a Banach space $\mathbb{E}$ of distributions that respects symmetries of the problem, we shall ask the norm of $\mathbb{E}$ to be invariant under translations in $\R^3$ ($\|f(x-x_0)\|_{\mathbb{E}}=\|f\|_\mathbb{E}$ for $x_0\in\R^3$) and under dilations ($\lambda \|f(\lambda x)\|_\mathbb{E}=\|f\|_\mathbb{E}$ for $\lambda>0$). In that case, Meyer \cite{Mey99} remarked that we have the embedding $\mathbb{E}\subset \dot B^{-1}_{\infty,\infty}$, a Besov space that plays a prominent role in the study of the Navier--Stokes equations. 

The use of Besov spaces in fluid mechanics relies essentially on the dyadic Littlewood--Paley decomposition (sometimes called the Littlewood--Paley--Stein  decomposition, as Stein is one of the first analysts to use it). This decomposition makes easy dealing with the non-linearity $\vu\cdot\vN\vu$ of the equations, by using the paraproduct operators of Bony \cite{Bon81}. Seminal works on Besov spaces and fluid mechanics appeared in the 90's, as the paper of Chemin in 1992 \cite{Che92} or the book of Cannone in 1995 \cite{Can95}; applications of the Littlewood-Paley decomposition to the borderline cases of regularity for solutions of Euler equations were given by Vishik in 1998-99 \cite{Vis98,Vis99}. Chemin developed a theory of time-space Besov spaces where the nonlinear evolution partial
differential equations are treated more efficiently after localization by
means of Littlewood–Paley decomposition \cite{ChL95, BCD11} (especially in the borderline cases of regularity).

An interesting example of the use of Besov spaces for Navier--Stokes equations is the proof of uniqueness of solutions in $\mathcal{C}([0,T],L^3)$  to the Cauchy problem for initial value $\vu^0\in L^3$. The first proof of such uniqueness has been given by Furioli, Lemari\'e-Rieusset and Terraneo \cite{FLT97, FLT00}.     If $\vu$ is a solution in $\mathcal{C}[0,T),L^3)$ and $\vu_K$ is the solution provided by Kato in $\mathcal{C}[0,T),L^3)$ with the additional property that $\lim_{t\rightarrow 0} t^{1/4} \|\vu_K\|_6=0$, then the function $\vw=\vu-\vu_K$ satisfies the identity
$$ \vw=-B(\vu_K,\vw)-B(\vw,\vu_K)-B(\vw,\vw).$$
The main step  of the proof of uniqueness  by  Furioli, Lemari\'e-Rieusset and Terraneo was to establish a contractive estimate (locally in time) on $\vw$ to prove that $\vw$ is equal to $0$, in spite of the fact that $B$ is not  bounded on $L^\infty L^3$.  They remarked that $\vw$ is more regular than $\vu$ and $\vu_K$ : $\vu-e^{t\Delta}\vu^0$ and  $\vu_K-e^{t\Delta}\vu^0$ belong to $L^\infty \dot B^{1/2}_{2,\infty}$; the contractive estimate  they found is then the following one : for $0<S<T$, 
$$ \sup_{0<t<S} \|\vw(t,.)\|_{\dot B^{1/2}_{2,\infty} }\leq   C  \sup_{0<t<S} \|\vw(t,.)\|_{\dot B^{1/2}_{2,\infty} } \left( \sup_{0<s<S} s^{1/8} \|\vu_K(s,.)\|_4+  \sup_{0<s<S} \|\vw(t,.)\|_{3} \right).$$ 
By continuity of both $\vu$ and $\vu_K$ in $L^3$ norm,   we have  $\lim_{t\rightarrow 0} \|\vw(t,.)\|_{L^{3,\infty}}=0$ and we find that we have a contractive estimate for $\vw$ if $S$ is small enough.  Hence, the fixed point $\vw$ is equal to $0$, and $\vu=\vu_K$ on $(0,S)$. The end of the proof follows by a bootstrap argument\footnote{This was after this result that Brezis asked me to write a book on Besov estimates for Navier-Stokes equations \cite{Lem02}.}.

In some points of the study of the Navier--Stokes equations, Besov spaces appear to be optimal. Let us quote three examples concerning the Leray solutions. We consider a solution $\vu\in L^\infty((0,T),L^2)\cap L^2((0,T),\dot H^1)$ of the Navier--Stokes equations with initial value $\vu^0\in L^2$, satisfying Leray's energy inequality :
$$ \|\vu(t,.)\|_2^2 + 2 \int_0^t \|\vN\otimes\vu\|_2^2\, ds\leq \|\vu^0\|_2^2.$$
\begin{itemize}
\item[$\bullet$]  regularity : a well-known result of Serrin \cite{Ser62} states that, if $\vu^0$ belongs more precisely to $H^1$, then $\vu$ will remain in $H^1$ as long as $\int_0^T \|\vu\|_q^p\, dt<+\infty$ with $3<q\leq\infty$ and $\frac 2 p+\frac 3 q=1$. The space $L^q$ has been replaced by many larger spaces with the same scaling properties. The largest one is $\dot B^{-\frac 3 q}_{\infty,\infty}$. Serrin's criterion has been proved to hold for $\vu\in L^p \dot B^{\sigma}_{\infty,\infty}$ for $\frac 2 p=1+\sigma$ and $1\leq p<+\infty$ (Kozono and Shimada \cite{KzS04} for $p>2$, Chen and Zhang \cite{ChZ06} for $1<p\leq 2$, Kozono, Ogawa and Taniuchi \cite{KOT02} for $p=1$).
\item[$\bullet$] weak-strong uniqueness : a well-known result of Prodi \cite{Pro59} and Serrin  \cite{Ser63} states that, if the Cauchy problem for $\vu^0$ has another solution $\vv$ in $L^\infty L^2\cap L^2 \dot H^1$ and if moreover $\vv\in L^p_t L^q_x$ with $\frac 2 p+\frac 3 q=1$ and $3<q\leq +\infty$, then $\vu=\vv$. Again, this has been extended by replacing $L^q$ by many larger spaces with the same scaling properties. Serrin's criterion has been proved to hold for $\vu\in L^p  X_\sigma$ for $\frac 2 p=1+\sigma$ and $1< p<+\infty$, where $X_\sigma=\dot B^\sigma_{\infty,\infty}$ if $\sigma>0$ (i.e. $p<2$) (Chen, Miao and Zhang \cite{CMZ09}), $X_0={\rm BMO}$ (Kozono and Taniuchi  \cite{KzT00}) and $X_\sigma=\dot M^{2,q}$ if $\sigma<0$ and $\sigma=-\frac 3 q$ (see the subsection on Morrey spaces).
\item[$\bullet$] energy (in)equality : a classical result of Lions \cite{Lio60} states that if the Leray solution $\vu$ satisfies $\vu\in L^4 L^4$, then Leray's energy inequality for $\vu$ is indeed an equality. The assumption $\vu\in L^4L^4$ has been weakened by Duchon and Robert \cite{DuR99} to $\vu\in L^3  b^{1/3}_{3,\infty}$ where $b^{1/3}_{3,\infty}$ is the closure of test functions in $\dot B^{1/3}_{3,\infty}$. (Remark that $L^2\dot H^1\cap L^4 L^4\subset L^3  b^{1/3}_{3,\infty}$).
\end{itemize}

\subsection*{Morrey spaces and Morrey--Campanato spaces.}
When dealing with scaled estimates in spaces of measurable functions, one is naturally driven to use Morrey spaces. The Morrey space $\dot M^{p,q}$, $1<p<+\infty$, $p\leq q\leq \infty$ is defined by :
$$ f\in \dot M^{p,q} \Leftrightarrow f\in L^p_{\rm loc} \text{ and }  \sup_{x_0\in\R^3, r>0} \frac 1{ \vert B(x_0,r)\vert^{\frac 1 p-\frac 1 q}} \| \1_{B(x_0,r)} f\|_p<+\infty.$$
Again, we define $\sigma=-\frac 3 q$  and we find equivalently
$$ f\in \dot M^{p,q} \Leftrightarrow  f\in L^p_{\rm loc} \text{ and }  \sup_{x_0\in\R^3, r>0} \frac 1{ r^{\frac 3 p+\sigma}} \| \1_{B(x,r)} f\|_p<+\infty.$$
The restriction $p\leq q\leq +\infty$ implies that we have   $-\frac 3 p\leq\sigma\leq 0$. Remark that if $\sigma<-\frac 3 p$ or $\sigma>0$ then $f=0$. Moreover, $\dot M^{p,\infty}=L^\infty$.

Morrey--Campanato spaces are quite similar, except that we correct $f$ with its mean value $m_{B(x_0,r)}f=\frac 1{\vert B(x_0,r)} \int_{B(x_0,r)} f\, dx$ :
$$ f\in  \mathcal{M}^{p,\sigma} \Leftrightarrow  f\in L^p_{\rm loc} \text{ and }  \sup_{x_0\in\R^3, r>0} \frac 1{ r^{\frac 3 p+\sigma}} \| \1_{B(x_0,r)} f -m_{B(x_0,r)}f\|_p<+\infty.$$
This time, $\sigma$ will be in the range $-\frac 3 p\leq \sigma<1$. Moreover, $  \mathcal{M}^{p,0}={\rm BMO}$ and, for $0<\sigma<1$, $  \mathcal{M}^{p,\sigma}=$ (the realization of) $\dot B^\sigma_{\infty,\infty}$ \cite{Cam63}. Thus, $\mathcal{M}^{p,\sigma}$ is the dual of the Hardy space $\mathcal{H}^r$ for $\frac 3 4<r\leq 1$ and  $\sigma=\frac 3 r-1$ \cite{FoS82}. 

If $\sigma<0$ and if $(\psi_{\epsilon,j,k})_{1\leq \epsilon\leq 7, j\in\Z,k\in Z^3}$ is a compactly supported wavelet bases with regularity  $\mathcal{C}^3$, then we find that 
$$ \vert \langle f\vert \psi_{\epsilon,j,k}\rangle \vert\leq   C 2^{3j(\frac 1 p-\frac 1 2)} \|\psi_\epsilon\|_{\frac p{p-1}}  \|f\|_{\mathcal{M}^{p,\sigma}} 2^{-j(\frac 3 p+\sigma)}$$
and we find that $\sum_{\epsilon,j,k}  \langle f\vert \psi_{\epsilon,j,k}\rangle \psi_{\epsilon,j,k}$ is (*-weakly) convergent in (the realization of)  $\dot B^\sigma_{\infty,\infty}$. Moreover, the series  $\sum_{\epsilon,j,k}  \langle f\vert \psi_{\epsilon,j,k}\rangle \vN \psi_{\epsilon,j,k}$  converges in $\mathcal{D}'$ to $\vN f$.  Thus, we have a decomposition $\mathcal{M}^{p,\sigma}=\dot M^{p,-\frac 3 \sigma}\oplus \R \1$ and an identification $\dot M^{p,-\frac 3 \sigma}= \mathcal{M}^{p,\sigma} \cap \dot B^\sigma_{\infty,\infty}$.

The first occurence of Morrey spaces in the study of Navier--Stokes equations was in a paper by
Giga and Miyakawa \cite{GiM89} on self--similar solutions. Then, in the early 90's, there has been results on mild solutions in Morrey spaces given by Kato \cite{Kat92}, Taylor 
\cite{Tay92} and
Federbush  \cite{Fed93}. In 1994, Kozono and Yamazaki \cite{KzY94} introduced Besov--Morrey spaces in order to give examples of singular initial values (or of initial values with large $L^3$ norms) leading to global mild solutions. Cannone's book \cite{Can95} or Lemari\'e-Rieusset's one \cite{Lem02} gave a systematic treatment of those spaces.

The flourishing of various classes of mild solutions for the Navier--Stokes equations that occured in the 90's opened the question of the largest space that would lead, through Picard iterations, to solutions. This space is included in $B^{-1}_{\infty,\infty}$ but is smaller, as the regularization by the heat kernel is not sufficient to give a meaning to the non-linear term. This space was identified by Koch and Tataru \cite{KoT01} and named ${\rm bmo}^{-1}$ : this is the space of distributions that are a sum of a bounded function $f_0\in L^\infty$ and of derivatives $\partial_j f_j$ of functions $f_j$ in the   ${\rm bmo}$ space of Goldberg (a local version of ${\rm BMO}$) \cite{Gol79}. The homogeneous version of this space is ${\rm BMO}^{-1}=\sqrt{-\Delta} ({\rm BMO})$. Recently, Auscher and Frey  \cite{AuF17} gave a new proof of the theorem of Koch and Tataru, based on the duality between the Hardy space $\mathcal{H}^1$ and ${\rm BMO}$. 

 Variations on the Koch and Tataru theorem led May \cite{May02} and Xiao \cite{Xia07} to consider initial values in $ (\sqrt{-\Delta})^{1-\sigma} \mathcal{M}^{2,-\sigma}= (\sqrt{-\Delta})^{1-\sigma} \dot M^{2,\frac 3 \sigma}$. Xiao linked his results to his theory of $Q$-spaces, and to Carleson measures and the tent spaces of Coifman, Meyer and Stein \cite{CMS85}.

Morrey spaces appear in many papers on the Navier--Stokes equations, extending results involving Lebesgue spaces where scaling properties prevail over global integrability. For instance, uniqueness of mild solutions in $\mathcal{C}([0,T],L^3)$ proven by Furioli, Lemari\'e--Rieusset and Terraneo holds as well in $\mathcal{C}([0,T], \dot m^{p,3})$ for $2<p\leq 3$, where $\dot m^{p,3}$ is the closure of test functions in the space $\dot M^{p,3}$ \cite{FLT00, Lem07}.{ The case of $\mathcal{C}([0,T], \dot m^{2,3})$ remains open. The space $\dot X^1=\mathcal{M}(\dot H^1\mapsto L^2)$ of pointwise multipliers from the Sobolev space $\dot H^1(\R^3)$ to $L^2(\R^3)$ satisfies the embeddings, for $2<p\leq 3$, $\dot M^{p,3}\subset \dot X^1\subset \dot M^{2,3}$ (Fefferman  \cite{Fef83}). May \cite{LMa07}  proved uniqueness of solutions in $\mathcal{C}([0,T], \dot  x^1)$, where $\dot x^1$  is the closure of test functions in the space $\dot  X^1$.

 Another interesting occurence of Morrey spaces in the study of the Navier--Stokes equations is the extension of the criterion of weak-strong uniqueness   of Prodi \cite{Pro59} and Serrin  \cite{Ser63}. The key point in the proof of the criterion is an inequality of the type
 $$ \vert \int  \vu\cdot (\vv\cdot\vN)\vv\, dx\vert\leq C \|\vu\|_{X_\sigma} \|\vv\|_2^{  {1+\sigma}} \|\vN\otimes\vv\|_2^{ {1-\sigma}}$$ for two divergence-free vector fields  $\vu$ and $\vv$. For $-1\leq \sigma\leq 0$, a simple approach is to use an inequality of the type $\| \vu\otimes\vv\|_2\leq C \|\vu\|_{X_\sigma} \|\vv\|_{H^{-\sigma}}$ together with $\|\vv\|_{H^{-\sigma}}\leq \|\vv\|_2^{1+\sigma} \|\vN\otimes\vv\|_2^{-\sigma}$. Let us write $\dot X^r=\mathcal{M}(\dot H^r\mapsto L^2)$ for the set of  pointwise multipliers from the Sobolev space $\dot H^r(\R^3)$ to $L^2(\R^3)$ (for a characterization of $\dot X^1$, see Maz'ya \cite{Maz64}); we get 
  $$ \vert \int  \vu\cdot (\vv\cdot\vN)\vv\, dx\vert\leq C \|\vu\|_{\dot X^r} \|\vv\|_2^{  {1-r}} \|\vN\otimes\vv\|_2^{ {1+r}}$$  and find weak-strong uniqueness for Leray solutions of the Cauchy problem with initial value $\vu_0$ if one of those solutions  belongs moreover to $L^p\dot X^r$ ($0\leq r<1$ and $\frac 2 p=1-r$) or to $\mathcal{C}([0,T), \dot x^1)$ (for $r=1$) \cite{Lem02} .

  For $0<r<1$, a better approach is to  use an inequality of the type $\| \vu\otimes\vv\|_2\leq C \|\vu\|_{X_\sigma} \|\vv\|_{\dot B^r_{2,1}}$ together with $\|\vv\|_{\dot B^r_{2,1}}\leq  C \|\vv\|_2^{1+r} \|\vN\otimes\vv\|_2^{r}$.  Thus, we are interested in the space $\mathcal{M}(\dot B^r_{2,1}\mapsto L^2)$   of  pointwise multipliers from the Besov space $\dot  B^r_{2,1}(\R^3)$ to $L^2(\R^3)$; this space turns out to be the Morrey space $\dot M^{2,\frac 3 r}$ \cite{Lem07, Lem16}, which is larger than $\dot X^r$. Hence, we get 
  $$ \vert \int  \vu\cdot (\vv\cdot\vN)\vv\, dx\vert\leq C \|\vu\|_{\dot M^{2,\frac 3 r}} \|\vv\|_2^{  {1-r}} \|\vN\otimes\vv\|_2^{ {1+r}}$$  and find weak-strong uniqueness for Leray solutions of the Cauchy problem with initial value $\vu_0$ if one of those solutions  belongs moreover to $L^p\dot  M^{2,\frac 3 r}$ ($0<  r<1$ and $\frac 2 p=1-r$). 
  
  For $\sigma\leq 0$, one uses the fact that $\vv$ is divergence free. 
  Recall that for $\sigma=0$, Kozono and Taniuchi  \cite{KzT00}wrote 
    $$ \vert \int  \vu\cdot (\vv\cdot\vN)\vv\, dx\vert\leq C \|\vu\|_{\rm BMO} \|\vv\cdot\vN\vv\|_{\mathcal{H}^1} \leq C'  \|\vu\|_{\rm BMO} \|\vv\|_2  \|\vN\otimes\vv\|_2$$ and  got 
   weak-strong uniqueness for Leray solutions of the Cauchy problem with initial value $\vu_0$ if one of those solutions  belongs moreover to  $L^2\, \rm BMO$.
  
  For $0<\sigma<1$, we use product laws in Sobolev spaces to estimate the (positive) regularity of $\vv\otimes\vv$ :
  $$ \|\vv\otimes\vv\|_{\dot B^{1-\sigma}_{1,1}} \leq C \|\vv\|_{\dot H^{\frac{1-\sigma}2}}^2 \leq C \|\vv\|_2^{1+\sigma} \|\vN\otimes \vv\|_2^{1-\sigma}$$ 
  so that
   $$  \vert \int  \vu\cdot (\vv\cdot\vN)\vv\, dx\vert\leq C \|\vN\otimes \vu\|_{\dot B^{\sigma-1}_{\infty,\infty}}  \|\vv\otimes\vv\|_{\dot B^{1-\sigma}_{1,1}} \leq C \|\vu\|_{\dot B^{\sigma}_{\infty,\infty}} \|\vv\|_2^{1+\sigma} \|\vN\otimes \vv\|_2^{1-\sigma}.$$  Thus,  find weak-strong uniqueness for Leray solutions of the Cauchy problem with initial value $\vu_0$ if one of those solutions  belongs moreover to $L^p\dot  B^\sigma_{\infty,\infty}$ ($0<  \sigma<1$ and $\frac 2 p=1+\sigma$).  The limit case $\sigma=1$ gives  weak-strong uniqueness when one of the solutions  belongs moreover to $L^1\, {\rm Lip}$.
  
 For $-1<\sigma<1$, those results may be unified in the following way :   weak-strong uniqueness for Leray solutions of the Cauchy problem with initial value $\vu_0$ if one of those solutions  belongs moreover to $L^p\mathcal{M}^{2,\sigma}$ ($-1<  \sigma<1$ and $\frac 2 p=1+\sigma$).  \\

Another point where scaling plays an important role is the theory of partial regularity for suitable  weak solutions of the Navier--Stokes solutions developed by Caffarelli, Kohn and Nirenberg \cite{CKN82}.
In order to simplify the proof given by Caffarelli, Kohn and Nirenberg in 1982, Ladyzhensakaya and Seregin  \cite{LaS99} used Morrey spaces as a basic tool for elliptic or parabolic equations. A systematic and inspiring proof wholly given ion terms of Morrey spaces has been given by Kukavica in 2011  \cite{Kuk11}.

 
 \section{Conclusion.}

 We have given many examples of   the interaction of harmonic analysis with the study of Navier--Stokes equations, beyond the simple use of Littlewood--Paley decomposition. The usefulness of such tools can be nicely illustrated by the case of the refined Gagliardo--Nirenberg inequalities of G\'erard, Meyer and Oru \cite{GMO96}. This inequality states that, if $1<p\leq +\infty$, $\alpha>0$ and $\beta>0$ then the control of $(\sqrt{-\Delta})^{\alpha}f$ in $L^p$ norm and of $f$ in $\dot B^{-\beta}_{\infty,\infty}$ gives a control of $f$ in $L^q$, with $\frac 1 q = \frac \beta {\alpha+\beta} \frac 1 p$ :
\begin{equation} \label{hedb} \|f\|_q\leq \|( \sqrt{-\Delta})^{\alpha}f\|_p^{\frac \beta  {\alpha+\beta}} 
 \|f\|_{\dot B^{-\beta}_{\infty,\infty}}^{\frac \alpha{\alpha+\beta}}\end{equation}
 The original proof is given in terms of the Littlewood--Paley decomposition of $f$ and of the characterization of $L^q$ as a Triebel--Lizorkin space $\dot F^0_{p,2}$.  But there is a very shorter and simpler proof based on Hedberg's inequality \cite{Hed72}. (More precisely a variant of Hedberg's inequality, where one replaces the role played by the Hardy--Littlewood maximal function  by Stein's maximal function, in order to be able to deal with distributions in $\dot B^{-\beta}_{\infty,\infty}$).  More precisely, if $N> \alpha/2$, one writes (for $f\in \dot B^{-\beta}_{\infty,\infty}$)
 $$ f= \frac{(-1)^N}{ \Gamma(N)}  \int_0^{+\infty} (t\Delta)^N e^{t\Delta} f\frac {dt}t$$
 and uses the inequalities 
 $$ \vert (t\Delta)^N e^{t\Delta} f(x)\vert\leq  C t^{\frac \alpha 2}\mathcal{M}_{(\sqrt{-\Delta})^\alpha f}(x) $$
 and 
  $$ \vert (t\Delta)^N e^{t\Delta} f(x)\vert\leq C t^{-\frac \beta 2} \|f\|_{\dot B^{-\beta}_{\infty,\infty}}$$
 to find Hedberg's inequality
 $$ \vert f(x)\vert \leq C \left(\mathcal{M}_{(\sqrt{-\Delta})^\alpha f}(x)\right)^{\frac \beta{\alpha+\beta}} \left( \|f\|_{\dot B^{-\beta}_{\infty,\infty}}\right)^{\frac \alpha{\alpha+\beta}}.$$
 Inequality (\ref{hedb}) is then obvious.
  
  Hedberg's inequality, combined with basic theory of singular integrals and maximal functions, should be a powerful tool to deal with some non-linear PDEs, avoiding  the rigidity of the Littlewood–Paley decomposition or of wavelet decompositions and in a way
replacing it by a molecular approach (molecules in Hardy spaces [where only size of the molecules is controlled) or in Besov spaces (where size and regularity of the molecules are controlled]). This was the claim in
 \cite{Lem12} and the basis for the book \cite{Lem16}.  Indeed, a Littlewood-Paley decomposition is stable neither through a transport
equation nor under the action of a singular integral convolution operator. On the other hand, a molecular decompostion will
be stable, since a molecule is preserved under a transport equation  with Lipschitzian drift (moving the center along the characteristic
curve and deforming the profile of the molecule, but without altering too much its scale), or through the action of
a singular integral convolution operator (with roughly speaking the same center and the same scale, but with a deformation of
the profile). Similarly, a wavelet decomposition is not preserved, but transformed into a vaguelette decomposition \cite{Lem02}. An interesting example of what can be done with molecules is the paper by Chamorro and Menozzi establishing regularization properties for an advection-diffusion problem with non-local diffusion and rough drift \cite{ChM18}; the title of their paper is quite programmatic for the use of real  methods in harmonic analysis when studying non-linear PDEs :  \textit{Non linear singular drifts and fractional operators: when Besov meets Morrey and Campanato},



\end{document}